\documentclass[11pt]{article}

\usepackage{amsfonts}
\usepackage{amssymb}
\usepackage{amsmath}
\usepackage{amsthm}
\usepackage{ngerman}

\def\II{I\hspace{-0.1cm}I}

\def\negthickspace{\!\!\!}

\newcommand{\nicefrac}[2]
{\leavevmode \kern.1em\raise.5ex\hbox{\the\scriptfont0 #1}
             \kern-.1em/\kern-.15em\lower.25ex
             \hbox{\the\scriptfont0 #2}}

\newtheorem*{theorem}{Theorem}
\newtheorem*{proposition}{Proposition}
\newtheorem*{corollary}{Corollary}

\newtheorem*{definition}{Definition}
\newtheorem*{lemma}{Lemma}

\theoremstyle{definition}
\newtheorem*{remark}{Remark}
\newtheorem*{remarks}{Remarks}

\theoremstyle{definition}

\begin{document} 

\begin{center}
{\Large{\sc $\mu$-stability of $2$-immersions}}\\[0.2cm]
{\Large{\sc of prescribed mean curvature}}\\[0.2cm]
{\Large{\sc and flat normal bundle}}\\[0.2cm]
{\Large{\sc in Euclidean spaces of higher codimension}}\\[1cm]
{\large Steffen Fr\"ohlich}\\[0.4cm]
{\small\bf Abstract}\\[0.4cm]
\begin{minipage}[c][2.5cm][l]{12cm}
{\small We present three ways to establish generalized stability inequalities for various classes of $2$-immersions in Euclidean spaces of higher codimension.}
\end{minipage}
\end{center}
{\small MCS 2000: 35J60, 53A07, 53A10}
\section{Basics and main results}
\subsection{$2$-immersions in $\mathbb R^n$}
On the closed unit disc $B:=\Big\{(u,v)\in\mathbb R^2\,:\,u^2+v^2\le 1\Big\}\subset\mathbb R^2$ we consider immersions
\begin{equation}\label{1.1}
  X=X(u,v)=(x^1(u,v),\ldots,x^n(u,v))\in C^3(B,\mathbb R^n)
\end{equation}
with the property $\mbox{rank}\,\partial X=2$ in $B$ for the Jacobian $\partial X\in\mathbb R^{n\times 2}.$\\[1ex]
Let ${\mathfrak N}=\{N_1,\ldots,N_{n-2}\}$ be a $C^2$-differentiable orthonormal section (ONS) of the normal bundle of the surface such that
\begin{equation}\label{1.2}
  X_{u^i}\cdot N_\sigma^t=0,\quad
  N_\sigma\cdot N_\omega^t=\delta_{\sigma\omega}
  \quad\mbox{for}\ i=1,2,\ \sigma,\omega=1,\ldots,n-2
\end{equation}
($u^1\equiv u,$ $u^2\equiv v$) with Kronecker's symbol $\delta_{\sigma\omega}.$ The superscript $t$ means the transposed vector.\\[1ex]
Let $N\in\mbox{span}\,{\mathfrak N}$ be a unit normal vector. We denote by
\begin{equation}\label{1.3}
\begin{array}{l}
  \displaystyle
  I(X)=(g_{ij})_{i,j=1,2}=(X_{u^i}\cdot X_{u^j}^t)_{i,j=1,2}\,, \\[0.2cm]
  \displaystyle
  \II_N(X)=(L_{N,ij})_{i,j=1,2}=-(X_{u^i}\cdot N_{u^j}^t)_{i,j=1,2}\,,
\end{array}
\end{equation}
the first fundamental form and the field of the second fundamental forms of $X$ w.r.t. $N.$
\subsection{Mean and Gaussian curvature fields}
In terms of these fundamental forms one defines the mean and Gaussian curvature fields by
\begin{equation}\label{1.4}
\begin{array}{lll}
  H_N\negthickspace
  & = & \negthickspace\displaystyle
        \frac{L_{N,11}g_{22}-2L_{N,12}g_{12}+L_{N,22}g_{11}}{g_{11}g_{22}-g_{12}^2}
        \,=\,\frac{\kappa_{N,1}+\kappa_{N,2}}{2}\,, \\[0.4cm]
  K_N\negthickspace
  & = & \negthickspace\displaystyle
        \frac{L_{N,11}L_{N,22}-L_{N,12}^2}{g_{11}g_{22}-g_{12}^2}
        \,=\,\kappa_{N,1}\kappa_{N,2}
\end{array}
\end{equation}
with the principle curvatures $\kappa_{N,1}$ and $\kappa_{N,2}$ w.r.t. the unit normal $N.$
\begin{remark}
$X$ is called a minimal surface iff $H_N\equiv 0$ in $B$ w.r.t. all unit normals $N.$
\end{remark}
\noindent
For a chosen ONS ${\mathfrak N}=\{N_1,\ldots,N_{n-2}\}$ we define (let $H_\sigma=H_{N_\sigma}$ etc.)
\begin{equation}\label{1.5}
  H:=\sqrt[+]{\sum_{\sigma=1}^{n-2}H_\sigma^2}\,,\quad
  K:=\sum_{\sigma=1}^{n-2}K_\sigma\,.
\end{equation}
In fact, $H$ and $K$ do not depend on the special choice of ${\mathfrak N}:$ $H$ is the length of the invariant mean curvature vector
\begin{equation}\label{1.6}
  \widehat N:=\sum_{\sigma=1}^{n-2}H_\sigma N_\sigma\,,
\end{equation}
and $K$ is the non-trivial component of the Riemannian curvature tensor.
\subsection{Conformal parameters}
We mainly use a conformal parametrization with the properties (see \cite{Sauvigny_03})
\begin{equation}\label{1.7}
  X_u^2=W=X_v^2\,,\quad X_u\cdot X_v^t=0\quad\mbox{in}\ B.
\end{equation}
\subsection{$\mu$-stability}
This paper presents the reader three ways to prove generalized stability in the following sense:
\begin{definition}
The conformally parametrized immersion $X=X(u,v)$ is called $\mu$-stable with a real constant $\mu\in(0,+\infty)$ and a function $q\in C^2(B,\mathbb R)$ iff it holds
\begin{equation}\label{1.8}
  \int\hspace{-0.25cm}\int\limits_{\hspace{-0.3cm}B}|\nabla\varphi|^2\,dudv
  \ge\mu\int\hspace{-0.25cm}\int\limits_{\hspace{-0.3cm}B}(q-K)W\varphi^2\,dudv
  \quad\mbox{for all}\ \varphi\in C_0^\infty(B,\mathbb R)
\end{equation}
where $q-K\ge 0$ in $B.$
\end{definition}
\begin{remarks}\quad
\begin{itemize}
\item[1.]
The definition is abstracted from the theory of the second variation: A conformally parametrized minimal surface $X$ (as a critical point of the area functional) is called stable iff
\begin{equation}\label{1.9}
  \int\hspace{-0.25cm}\int\limits_{\hspace{-0.3cm}B}|\nabla\varphi|^2\,dudv
  \ge 2\int\hspace{-0.25cm}\int\limits_{\hspace{-0.3cm}B}(-K)W\varphi^2\,dudv.
\end{equation}
This means $\mu$-stability with $\mu=2$ and $q\equiv 0$ (note that $K\le 0$). And a conformally parametrized surface with constant mean curvature $h_0\in\mathbb R$ (as a critical point of the area functional with a suitable volume constraint) is stable iff
\begin{equation}\label{1.10}
  \int\hspace{-0.25cm}\int\limits_{\hspace{-0.3cm}B}|\nabla\varphi|^2\,dudv
  \ge 2\int\hspace{-0.25cm}\int\limits_{\hspace{-0.3cm}B}(2h_0^2-K)W\varphi^2\,dudv,
\end{equation}
that is, it is $\mu$-stable with $\mu=2$ and $q\equiv 2h_0^2.$
\item[2.]
Various criteria are known to realize $\mu$-stability in $\mathbb R^3:$ For example, stability for minimal surfaces \cite{Barbosa_doCarmo_01}, for surfaces of prescribed constant mean curvature \cite{Ruchert_01}, for F-minimal surfaces \cite{Clarenz_01}, or for weighted minimal surfaces \cite{Froehlich_01}.
\item[3.]
In \cite{Barbosa_doCarmo_02} the reader can find stability criteria for minimal surfaces in the three-sphere $S^3,$ in the hyperbolic space $H^3,$ and in the Euclidean space $R^n.$ In our paper we will especially present a new variation of this last result.
\end{itemize}
\end{remarks}
\subsection{Applications of $\mu$-stability}
\begin{itemize}
\item[1.]
Stability conditions are used to control the area growth of geodesic discs or Dirichlet energies of their Gauss mappings. In turn, these estimates help to establish moduli of continuity for the mappings, and, finally, such moduli are needed for inner Schauder estimates. We refer to \cite{Bergner_Froehlich_01}, \cite{Froehlich_01}, \cite{Froehlich_02}, \cite{Froehlich_03}, \cite{Sauvigny_01}, \cite{Sauvigny_02} for such applications.
\vspace*{-1.2ex}
\item[2.]
Stability conditions are the basis for curvature estimates following the methods of \cite{Ecker_Huisgen_01} and \cite{Schoen_Simon_Yau_01}; see \cite{Froehlich_Winklmann_01}, \cite{Winklmann_01}, \cite{Winklmann_02}, and also \cite{Clarenz_01}.
\end{itemize}
\subsection{Main results}
Let us enumerate the results of this paper:
\begin{itemize}
\item[$\to$]
{\it Section 2: }We compute the first and second variation of the Fermat type functional
\begin{equation}\label{1.11}
  {\mathcal F}[X]:=\int\hspace{-0.25cm}\int\limits_{\hspace{-0.3cm}B}\Gamma(X)W\,dudv
  \longrightarrow\mbox{extr!}
\end{equation}
with a weight $\Gamma\in C^1(\mathbb R^n,\mathbb R_+),$ given such that
\begin{equation}\label{1.12}
  0<\Gamma_{min}\le\Gamma(X)\le\Gamma_{max}<+\infty,
\end{equation}
and we derive $\mu$-stability criteria for stable critical points;
\vspace*{-0.4ex}
\item[$\to$]
{\it Section 3: }We establish $\mu$-stability for surface graphs of prescribed mean curvature fields;
\vspace*{-0.4ex}
\item[$\to$]
{\it Section 4: }Differential equations for the Hopf functions are derived which give informations about the zeros of the Gaussian curvature;
\vspace*{-0.4ex}
\item[$\to$]
{\it Section 5: }$\mu$-stability for minimal immersions with flat normal bundles is established following methods of Ruchert \cite{Ruchert_01} and Barbosa/do\,Carmo \cite{Barbosa_doCarmo_02}.
\end{itemize}
\section{The second variation and stability}
\setcounter{equation}{0}
Consider the unit normal field (summation over $\sigma$!)
\begin{equation}\label{2.1}
  N_{\widehat\gamma}:=\widehat\gamma^\sigma N_\sigma,\quad
  \sum_{\sigma=1}^{n-2}(\widehat\gamma^\sigma)^2=1
  \quad\mbox{in}\ B
\end{equation}
with coefficients $\widehat\gamma^\sigma\in C^2(B,\mathbb R),$ and the variation
\begin{equation}\label{2.2}
  \widetilde X=X+\varepsilon\varphi N_{\widehat\gamma}\quad\mbox{in}\ B,\quad
  \varphi\in C_0^\infty(B,\mathbb R),\ \varepsilon\in(-\varepsilon_0,+\varepsilon_0).
\end{equation}
\subsection{First and second variation of the area element}
\begin{lemma}
Using conformal parameters $(u,v)\in B,$ there hold
\begin{equation}\label{2.3}
  \delta_{N_{\widehat\gamma}}h_{11}=-2\varphi L_{\widehat\gamma,11}\,,\quad
  \delta_{N_{\widehat\gamma}}h_{12}=-2\varphi L_{\widehat\gamma,12}\,,\quad
  \delta_{N_{\widehat\gamma}}h_{22}=-2\varphi L_{\widehat\gamma,22}
\end{equation}
for the variation w.r.t. $N_{\widehat\gamma},$ as well as
\begin{equation}\label{2.4}
\begin{array}{lll}
  \delta_{N_{\widehat\gamma}}^2h_{11}\negthickspace
  & = & \negthickspace\displaystyle
        2\varphi_u^2
        +\frac{2}{W}\,\varphi^2\,\big\{L_{\widehat\gamma,11}^2+L_{\widehat\gamma,12}^2\big\}
        +2\varphi^2
          \sum_{\sigma=1}^{n-2}
          \Big\{\widehat\gamma_u^\sigma+\widehat\gamma^\omega T_{\omega,1}^\sigma\Big\}^2\,, \\[0.6cm]
  \delta_{N_{\widehat\gamma}}^2h_{12}\negthickspace
  & = & \negthickspace\displaystyle
        2\varphi_u\varphi_v
        +\frac{2}{W}\,\varphi^2(L_{\widehat\gamma,11}+L_{\widehat\gamma,22})L_{\widehat\gamma,12}
        +2\varphi^2
          \sum_{\sigma=1}^{n-2}
          \Big\{\widehat\gamma_u^\sigma+\widehat\gamma^\omega T_{\omega,1}^\sigma\Big\}
          \Big\{\widehat\gamma_v^\sigma+\widehat\gamma^\omega T_{\omega,2}^\sigma\Big\}, \\[0.6cm]
  \delta_{N_{\widehat\gamma}}^2h_{22}\negthickspace
  & = & \negthickspace\displaystyle
        2\varphi_v^2
        +\frac{2}{W}\,\varphi^2\,\big\{L_{\widehat\gamma,12}^2+L_{\widehat\gamma,22}^2\big\}
        +2\varphi^2
          \sum_{\sigma=1}^{n-2}
          \Big\{\widehat\gamma_v^\sigma+\widehat\gamma^\omega T_{\omega,2}^\sigma\Big\}^2
\end{array}
\end{equation}
with the torsion coefficients $T_{\sigma,k}^\omega$ defined in (\ref{2.9}).
\end{lemma}
\begin{proof}
\begin{itemize}
\item[1.]
There hold $\widetilde X_u=X_u+\varepsilon\varphi_uN_{\widehat\gamma}+\varepsilon\varphi N_{\widehat\gamma,u},$ $\widetilde X_v=X_v+\varepsilon\varphi_vN_{\widehat\gamma}+\varepsilon\varphi N_{\widehat\gamma,v},$ and, therefore,
\begin{equation}\label{2.5}
\begin{array}{rcl}
  \widetilde X_u^2\negthickspace
  & = & \negthickspace
        W+2\varepsilon\varphi X_u\cdot N_{\widehat\gamma,u}^t
        +\varepsilon^2\varphi_u^2
        +\varepsilon^2\varphi^2N_{\widehat\gamma,u}^2\,, \\[0.2cm]
  \widetilde X_v\negthickspace
  & = & \negthickspace
        W+2\varepsilon\varphi X_v\cdot N_{\widehat\gamma,v}^t
        +\varepsilon^2\varphi_v^2
        +\varepsilon^2\varphi^2N_{\widehat\gamma,v}^2\,, \\[0.2cm]
  \widetilde X_u\cdot\widetilde X_v^t\negthickspace
  & = & \negthickspace
        \varepsilon\big\{X_u\cdot N_{\widehat\gamma,v}^t+X_v\cdot N_{\widehat\gamma,u}^t\big\}\,\varphi
        +\varepsilon^2\varphi_u\varphi_v
        +\varepsilon^2\varphi^2N_{\widehat\gamma,u}\cdot N_{\widehat\gamma,v}^t
\end{array}
\end{equation}
taking $N_{\widehat\gamma,u^i}\cdot N_{\widehat\gamma}^t=0,$ $i=1,2,$ into account.
\item[2.]
We introduce the forms
\begin{equation}\label{2.6}
  L_{\widehat\gamma,ij}
  :=X_{u^iu^j}\cdot N_{\widehat\gamma}^t
   =-X_{u^i}\cdot N_{\widehat\gamma,u^j}^t
   =-X_{u^j}\cdot N_{\widehat\gamma,u^i}^t
\end{equation}
noting that $X_{u^i}\cdot N_{\widehat\gamma,u^j}^t=-X_{u^iu^j}\cdot N_{\widehat\gamma}^t,$ and $X_u\cdot N_{\widehat\gamma,v}^t=X_v\cdot N_{\widehat\gamma,u}^t$ due to the symmetry of $L_{\sigma,ij}.$ Then, equations (\ref{2.5}) can be written as
\begin{equation}\label{2.7}
\begin{array}{rcl}
  \widetilde X_u^2\negthickspace
  & = & \negthickspace
        W-2\varepsilon\varphi L_{\widehat\gamma,11}
        +\varepsilon^2\varphi_u^2+\varepsilon^2\varphi^2N_{\widehat\gamma,u}^2\,, \\[0.2cm]
  \widetilde X_v^2\negthickspace
  & = & \negthickspace
        W-2\varepsilon\varphi L_{\widehat\gamma,22}
        +\varepsilon^2\varphi_v^2+\varepsilon^2\varphi^2N_{\widehat\gamma,v}^2\,, \\[0.2cm]
  \widetilde X_u\cdot\widetilde X_v^t\negthickspace
  & = & \negthickspace
        -\,2\varepsilon\varphi L_{\widehat\gamma,12}
        +\varepsilon^2\varphi_u\varphi_v+\varepsilon^2\varphi^2N_{\widehat\gamma,u}\cdot N_{\widehat\gamma,v}^t\,.
\end{array}
\end{equation}
\item[3.]
The Weingarten equations (see \cite{Brauner_01})
\begin{equation}\label{2.8}
  N_{\sigma,u^i}^t=-L_{\sigma,ij}h^{jk}\,X_{u^k}^t+T_{\sigma,i}^\omega N_\omega^t\,,\quad
  i=1,2,\ \sigma=1,\ldots,n-2,
\end{equation}
with the torsion coefficients
\begin{equation}\label{2.9}
  T_{\sigma,i}^\omega
  :=\left\{
      \begin{array}{cl}
        N_{\sigma,u^i}\cdot N_\omega^t\,, & \mbox{if}\ \sigma\not=\omega \\[0.2cm]
        0, & \mbox{else}
      \end{array}
    \right.,
\end{equation}
and written in terms of conformal parameters, imply
\begin{equation}\label{2.10}
\begin{array}{lll}
  N_{\widehat\gamma,u}\negthickspace
  & = & \negthickspace\displaystyle
        -\,\frac{L_{\widehat\gamma,11}}{W}\,X_u
        -\frac{L_{\widehat\gamma,12}}{W}\,X_v
        +\Big\{\widehat\gamma_u^\sigma+\widehat\gamma^\omega T_{\omega,1}^\sigma\Big\}\,N_\sigma\,, \\[0.6cm]
  N_{\widehat\gamma,v}\negthickspace
  & = & \negthickspace\displaystyle
        -\,\frac{L_{\widehat\gamma,12}}{W}\,X_u
        -\frac{L_{\widehat\gamma,22}}{W}\,X_v
        +\Big\{\widehat\gamma_v^\sigma+\widehat\gamma^\omega T_{\omega,2}^\sigma\Big\}\,N_\sigma\,.
\end{array}
\end{equation}
This in turn leads to
\begin{equation}\label{2.11}
\begin{array}{rcl}
  N_{\widehat\gamma,u}^2\negthickspace
  & = & \negthickspace\displaystyle
        \frac{L_{\widehat\gamma,11}^2+L_{\widehat\gamma,12}^2}{W}
        +\sum_{\sigma=1}^{n-2}
         \Big\{\widehat\gamma_u^\sigma+\widehat\gamma^\omega T_{\omega,1}^\sigma\Big\}^2\,, \\[0.6cm]
  N_{\widehat\gamma,u}\cdot N_{\widehat\gamma,v}^t\negthickspace
  & = & \negthickspace\displaystyle
        \frac{L_{\widehat\gamma,11}+L_{\widehat\gamma,22}}{W}\,L_{\widehat\gamma,12}
        +\sum_{\sigma=1}^{n-2}
         \Big\{\widehat\gamma_u^\sigma+\widehat\gamma^\omega T_{\omega,1}^\sigma\Big\}
         \Big\{\widehat\gamma_v^\sigma+\widehat\gamma^\omega T_{\omega,2}^\sigma\Big\}, \\[0.6cm]
  N_{\widehat\gamma,v}^2\negthickspace
  & = & \negthickspace\displaystyle
        \frac{L_{\widehat\gamma,12}^2+L_{\widehat\gamma,22}^2}{W}
        +\sum_{\sigma=1}^{n-2}
         \Big\{\widehat\gamma_v^\sigma+\widehat\gamma^\omega T_{\omega,2}^\sigma\Big\}^2\,.
\end{array}
\end{equation}
Inserting into (\ref{2.7}) proves the statement.
\end{itemize}
\end{proof}
\noindent
Note that $2W\delta W=h_{22}\delta h_{11}-2h_{12}\delta h_{12}+h_{11}\delta h_{22}.$ Thus, $\delta W=-2H_{\widehat\gamma}W\varphi$ in view of (\ref{2.3}) with the mean curvature $H_{\widehat\gamma}$ w.r.t. $N_{\widehat\gamma}.$ Together with
\begin{equation}\label{2.12}
  \delta^2W
  =\frac{1}{2}\,\delta^2h_{11}+\frac{1}{2}\,\delta^2h_{22}
   +\frac{1}{W}\,\delta h_{11}\delta h_{22}
   -\frac{1}{W}\,(\delta h_{12})^2
   -\frac{1}{W}\,(\delta W)^2
\end{equation}
using conformal parameters we arrive at the
\begin{corollary}
(Second variation of the area element)\\
Using conformal parameter it holds
\begin{equation}\label{2.13}
  \delta^2W
  =|\nabla\varphi|^2
   +2K_{\widehat\gamma}W\varphi^2
   +\sum_{\sigma=1}^{n-2}
         \left\{
           \Big(\widehat\gamma_u^\sigma+\widehat\gamma^\omega T_{\omega,1}^\sigma\Big)^2
           +\Big(\widehat\gamma_v^\sigma+\widehat\gamma^\omega T_{\omega,2}^\sigma\Big)^2
         \right\}\varphi^2
\end{equation}
for arbitrary $\varphi\in C_0^\infty(B,\mathbb R)$ with the Gaussian curvature $K_{\widehat\gamma}$ w.r.t. $N_{\widehat\gamma}.$
\end{corollary}
\subsection{First and second variation of Fermat's functional}
We compute $\delta[\Gamma(X)W]=[\delta\Gamma(X)]W+\Gamma(X)\,\delta W=\Gamma_X(X)\cdot N_{\widehat\gamma}^t\,W\varphi-2\Gamma(X)H_{\widehat\gamma}W\varphi\,.$ Thus, for a critical point of Fermat's functional it holds
\begin{equation}\label{2.14}
  0=\delta\int\hspace{-0.25cm}\int\limits_{\hspace{-0.3cm}B}\Gamma(X)W\,dudv
   =\int\hspace{-0.25cm}\int\limits_{\hspace{-0.3cm}B}
    \Big\{
      \Gamma_X(X)\cdot N_{\widehat\gamma}^t-2\Gamma(X)H_{\widehat\gamma}
    \Big\}\,W\varphi\,dudv
\end{equation}
for all $\varphi\in C_0^\infty(B,\mathbb R),$ which gives its mean curvature w.r.t. $N_{\widehat\gamma}$
\begin{equation}\label{2.15}
  H_{\widehat\gamma}=\frac{\Gamma_X(X)\cdot N_{\widehat\gamma}^t}{2\Gamma(X)}\,.
\end{equation}
\begin{theorem}
The second variation of Fermat's functional ${\mathcal F}[X]$ is
\begin{equation}\label{2.16}
\begin{array}{lll}
  \displaystyle
  \delta_{N_{\widehat\gamma}}^2{\mathcal F}[X;\varphi]\negthickspace
  & = & \negthickspace\displaystyle
        \int\hspace{-0.25cm}\int\limits_{\hspace{-0.3cm}B}
        \Gamma(X)|\nabla\varphi|^2\,dudv \\[0.7cm]
  &   & \negthickspace\displaystyle
        +\,2\int\hspace{-0.25cm}\int\limits_{\hspace{-0.3cm}B}
            \Big\{
              H_{\widehat\gamma,X}\cdot N_{\widehat\gamma}^t
              -2H_{\widehat\gamma}^2
              +K_{\widehat\gamma}
            \Big\}\,\Gamma(X)W\varphi^2\,dudv \\[0.7cm]
  &   & \negthickspace\displaystyle
        +\,\int\hspace{-0.25cm}\int\limits_{\hspace{-0.3cm}B}
           \sum_{\sigma=1}^{n-2}
           \left\{
             \Big(\widehat\gamma_u^\sigma+\widehat\gamma^\omega T_{\omega,1}^\sigma\Big)^2
             +\Big(\widehat\gamma_v^\sigma+\widehat\gamma^\omega T_{\omega,2}^\sigma\Big)^2
           \right\}\Gamma(X)\varphi^2\,dudv
\end{array}
\end{equation}
for arbitrary $\varphi\in C_0^\infty(B,\mathbb R)$ along $N_{\widehat\gamma}=\widehat\gamma^\sigma N_\sigma,$ and the mean curvature field $H_{\widehat\gamma}$ from (\ref{2.15}).
\end{theorem}
\begin{proof}
From $\frac{\partial}{\partial\varepsilon}\,\Gamma(\widetilde X)\widetilde W=\Gamma_X(\widetilde X)\cdot N_{\widehat\gamma}^t\,\widetilde W\varphi+\Gamma(\widetilde X)\,\frac{\partial}{\partial\varepsilon}\,\widetilde W$ we conclude
\begin{equation}\label{2.17}
  \frac{\partial}{\partial\varepsilon}\,{\mathcal F}[\widetilde X]
  =\int\hspace{-0.25cm}\int\limits_{\hspace{-0.3cm}B}
   \left\{
     2\Gamma(\widetilde X)H(\widetilde X,N_{\widehat\gamma})\widetilde W\varphi
     +\Gamma(\widetilde X)\,\frac{\partial}{\partial\varepsilon}\,\widetilde W
   \right\}\,dudv.
\end{equation}
A further differentiation w.r.t. $\varepsilon$ at $\varepsilon=0$ gives
\begin{equation}\label{2.18}
\begin{array}{lll}
  \displaystyle
  \delta_{N_{\widehat\gamma}}^2{\mathcal F}[X;\varphi]\negthickspace\displaystyle
  & = & \negthickspace\displaystyle
        2\int\hspace{-0.25cm}\int\limits_{\hspace{-0.3cm}B}
         \Big\{
           \Gamma_X(X)\cdot N_{\widehat\gamma}^t\,H(X,N_{\widehat\gamma})
           +\Gamma(X)H_X(X,N_{\widehat\gamma})\cdot N_{\widehat\gamma}^t
         \Big\}W\varphi^2\,dudv \\[0.8cm]
  &   & \negthickspace\displaystyle
        -\,2\int\hspace{-0.25cm}\int\limits_{\hspace{-0.3cm}B}
            \Big\{
              2\Gamma(X)H(X,N_{\widehat\gamma})^2
              +\Gamma_X(X)\cdot N_{\widehat\gamma}^t\,H(X,N_{\widehat\gamma})
            \Big\}W\varphi^2\,dudv \\[0.8cm]
  &   & \negthickspace\displaystyle
        +\int\hspace{-0.25cm}\int\limits_{\hspace{-0.3cm}B}\Gamma(X)|\nabla\varphi|^2\,dudv
        +2\int\hspace{-0.25cm}\int\limits_{\hspace{-0.3cm}B}\Gamma(X)K_{\widehat\gamma}W\varphi^2\,dudv \\[0.8cm]
  &   & \negthickspace\displaystyle
        +\int\hspace{-0.25cm}\int\limits_{\hspace{-0.3cm}B}
         \sum_{\sigma=1}^{n-2}
         \left\{
           \Big(\widehat\gamma_u^\sigma+\widehat\gamma^\omega T_{\omega,1}^\sigma\Big)^2
           +\Big(\widehat\gamma_v^\sigma+\widehat\gamma^\omega T_{\omega,2}^\sigma\Big)^2
         \right\}\Gamma(X)\varphi^2\,dudv,
\end{array}
\end{equation}
and the statement follows.
\end{proof}
\begin{remarks}
\begin{itemize}
\item[1.]
In case $n=3$ all torsion coefficients vanish, and it follows that
\begin{equation}\label{2.19}
  \delta_N^2{\mathcal F}[X;\varphi]
  =\int\hspace{-0.25cm}\int\limits_{\hspace{-0.3cm}B}
   \Big\{
     |\nabla\varphi|^2+2H_X(X,N)\cdot N^t-4H(X,N)^2+2K
   \Big\}\,\Gamma(X)\,dudv
\end{equation}
for all $\varphi\in C_0^\infty(B,\mathbb R).$ Now, let $X$ be stable, that is $\delta_N^2{\mathcal F}[X;\varphi]\ge 0$ for all $\varphi\in C_0(B,\mathbb R).$
\end{itemize}
\begin{corollary}
A stable critical point $X$ of Fermat's functional, such that $2H^2-H_X\cdot N^t-K\ge 0,$ is $\mu$-stable with $0<\mu\le\frac{2\Gamma_{min}}{\Gamma_{max}}$ and $q=2H^2-H_X\cdot N^t,$ $H(X,Z)$ from (\ref{2.15}).
\end{corollary}
\begin{itemize}
\item[2.]
Let $X$ be a minimal surface, that is $X$ is critical for the area functional ($\Gamma(X)\equiv 1$). Set $\widehat\gamma^\vartheta=1$ for one index $\vartheta\in\{1,\ldots,n-2\},$ and $\widehat\gamma^\sigma=0$ for $\sigma\not=\vartheta.$ Then
\begin{equation}\label{2.20}
  \delta_{N_\vartheta}^2{\mathcal A}[X]
  =\int\hspace{-0.25cm}\int\limits_{\hspace{-0.3cm}B}
   (|\nabla\varphi|^2+2K_\vartheta W\varphi^2)\,dudv
   +\int\hspace{-0.25cm}\int\limits_{\hspace{-0.3cm}B}
    \sum_{\sigma=1}^{n-2}
    \Big\{
      (T_{\vartheta,1}^\sigma)^2+(T_{\vartheta,2}^\sigma)^2
    \Big\}\,\varphi^2\,dudv.
\end{equation}
Assume $X$ is stable, that is $\delta_{N_\vartheta}^2{\mathcal A}[X]\ge 0$ for all $\vartheta=1,\ldots,n-2$ w.r.t. the an ONS ${\mathfrak N}.$ Then summing up these $n-2$ stability inequalities (with the same $\varphi$) gives
\begin{equation}\label{2.21}
\begin{array}{lll}
  \displaystyle
  \int\hspace{-0.25cm}\int\limits_{\hspace{-0.3cm}B}|\nabla\varphi|^2\,dudv\negthickspace
  & \ge & \negthickspace\displaystyle
          \frac{2}{n-2}
          \int\hspace{-0.25cm}\int\limits_{\hspace{-0.3cm}B}
          (-K)W\varphi^2\,dudv \\[0.8cm]
  &     & \negthickspace\displaystyle
          -\,\frac{1}{n-2}
             \int\hspace{-0.25cm}\int\limits_{\hspace{-0.3cm}B}
             \sum_{\sigma,\vartheta=1}^{n-2}
             \Big\{
               (T_{\vartheta,1}^\sigma)^2+(T_{\vartheta,2}^\sigma)^2
             \Big\}\,\varphi^2\,dudv
\end{array}
\end{equation}
for all $\varphi\in C_0^\infty(B,\mathbb R).$ The second integral is what we call the total torsion of ${\mathfrak N}.$ If $X$ has flat normal bundle, then there exists an ONS ${\mathfrak N}$ which is free of torsion: $T_{\vartheta,k}^\sigma\equiv 0.$
\end{itemize}
\begin{corollary}
The stable minimal surface $X\in{\mathfrak C}(B,\mathbb R^n)$ with flat normal bundle is $\mu$-stable with $0<\mu\le\frac{2}{n-2}$ and $q\equiv 0.$
\end{corollary}
\noindent
The general case of non-flat normal bundle will be analysed in a future paper.
\end{remarks}
\section{$\mu$-stability due to Schwarz}
\setcounter{equation}{0}
\subsection{A differential equation for the projection mapping $\chi$}
The following function $\chi$ is the initial point of Schwarz' method to prove stability.
\begin{lemma}
Using conformal parameters, the function
\begin{equation}\label{3.1}
  \chi:=\frac{J}{W}\,,\quad J:=x_u^1x_v^2-x_u^2x_v^1\,,
\end{equation}
satisfies in $B$ the differential equation
\begin{equation}\label{3.2}
\begin{array}{lll}
  \triangle\chi\negthickspace
  & = & \negthickspace\displaystyle
        -\,2\sum_{\sigma=1}^{n-2}
            \Big\{
              2H_\sigma^2-K_\sigma
            \Big\}\,W\chi
        +\sum_{\sigma=1}^{n-2}\sum_{\omega=1}^{n-2}
         S_{\sigma,12}^\omega(n_\sigma^1n_\omega^2-n_\sigma^2n_\omega^1) \\[0.5cm]
  &   & \negthickspace\displaystyle
        +\,2\sum_{\sigma=1}^{n-2}\sum_{\omega=1}^{n-2}
            \Big\{
              (n_\sigma^1x_v^2-n_\sigma^2x_v^1)T_{\sigma,1}^\omega
              -(n_\sigma^1x_u^2-n_\sigma^2x_u^1)T_{\sigma,2}^\omega
            \Big\}\,H_\sigma \\[0.5cm]
  &   & \negthickspace\displaystyle
        +\,2\sum_{\sigma=1}^{n-2}
            \Big\{
              H_{\sigma,u}(n_\sigma^1x_v^2-n_\sigma^2x_v^1)
              -H_{\sigma,v}(n_\sigma^1x_u^2-n_\sigma^2x_u^1)
            \Big\}
\end{array}
\end{equation}
with the components
\begin{equation}\label{3.3}
  S_{\sigma,ij}^\omega
  :=T_{\sigma,i,u^j}^\omega-T_{\sigma,j,u^i}^\omega
    +T_{\sigma,i}^\theta T_{\theta,j}^\omega
    -T_{\sigma,j}^\theta T_{\theta,i}^\omega\,,\quad
  i,j=1,2,\ \sigma,\omega=1,\ldots,n-2,
\end{equation}
of the curvature tensor of the normal bundle.
\end{lemma}
\begin{remarks}
\begin{itemize}
\item[1.]
For immersions in $\mathbb R^3,$ the function $\chi$ equals the third component of the unit normal vector of the surface. Thus, we call $\chi$ the projection mapping of the graph.
\item[2.]
Due to Ricci's integrability (see \cite{doCarmo_01}) conditions there hold 
\begin{equation}\label{3.4}
  S_{\sigma,12}^\omega=g^{jk}(L_{\sigma,1j}L_{\omega,k2}-L_{\sigma,2j}L_{\omega,k1})
  \quad\mbox{for}\ \sigma,\omega=1,\ldots,n-2.
\end{equation}
\end{itemize}
\end{remarks}
\def\proofname{Sketch of the proof}
\begin{proof}\quad
\begin{itemize}
\item[1.]
We compute the first derivatives $\chi_u$ and $\chi_v:$ Using the Gauss equations (see \cite{Brauner_01})
\begin{equation}\label{3.5}
  X_{u^iu^j}=\Gamma_{ij}^k\,X_{u^k}+\sum_{\sigma=1}^{n-2}L_{\sigma,ij}N_\sigma\,,\quad i,j=1,2,
\end{equation}
with the Christoffel symbols $\Gamma_{ij}^k=\frac{1}{2}\,g^{k\ell}(g_{\ell i,u^j}+g_{j\ell,u^i}-g_{ij,u^\ell}),$ $i,j,k=1,2,$ we calculate
\begin{equation}\label{3.6}
\begin{array}{lll}
  \chi_u\negthickspace
  & = & \negthickspace\displaystyle
        \frac{1}{W}\,(x_{uu}^1x_v^2+x_u^1x_{uv}^2-x_{uu}^2x_v^1-x_u^2x_{uv}^1)
        -\frac{W_u}{W^2}\,J \\[0.4cm]
  & = & \negthickspace\displaystyle
        \frac{1}{2W^2}\,
        (W_ux_u^1x_v^2-W_vx_v^1x_v^2+W_vx_u^2x_u^1+W_ux_v^2x_u^1) \\[0.5cm]
  &   & \negthickspace\displaystyle
        +\,\frac{1}{2W^2}\,(-W_ux_u^2x_v^1+W_vx_v^2x_v^1-W_vx_u^1x_u^2-W_ux_v^1x_u^2)
        -\frac{W_u}{W^2}\,J \\[0.4cm]
  &   & \negthickspace\displaystyle
        +\,\frac{1}{W}\,
           \sum_{\sigma=1}^{n-2}
           \Big\{
             L_{\sigma,11}(n_\sigma^1x_v^2-n_\sigma^2x_v^1)
             -L_{\sigma,12}(n_\sigma^1x_u^2-n_\sigma^2x_u^1)
           \Big\} \\[0.5cm]
  & = & \negthickspace\displaystyle
        \frac{1}{W}\,
        \sum_{\sigma=1}^{n-2}
           \Big\{
            L_{\sigma,11}(n_\sigma^1x_v^2-n_\sigma^2x_v^1)
             -L_{\sigma,12}(n_\sigma^1x_u^2-n_\sigma^2x_u^1)
           \Big\}
\end{array}
\end{equation}
using conformal parameters, and analogously for $\chi_v.$
\item[2.]
A computation of $\chi_{uu},$ taking $-L_{\sigma,22,u}+2H_{\sigma,u}W+2H_\sigma W_u,$ the Weingarten equations and the Gauss equations into account, gives
\begin{equation}\label{3.7}
\begin{array}{lll}
  \chi_{uu}\negthickspace
  & = & \negthickspace\displaystyle
        -\,\frac{W_u}{W^2}\,
           \sum_{\sigma=1}^{n-2}
           \Big\{
             L_{\sigma,11}(n_\sigma^1x_v^2-n_\sigma^2x_v^1)
             -L_{\sigma,12}(n_\sigma^1x_u^2-n_\sigma^2x_u^1)
           \Big\} \\[0.5cm]
  &   & \negthickspace\displaystyle
        -\,\frac{1}{W}\,
           \sum_{\sigma=1}^{n-2}
           \Big\{
             L_{\sigma,22,u}(n_{\sigma}^1x_v^2-n_{\sigma}^2x_v^1)
             +L_{\sigma,12,u}(n_{\sigma}^1x_u^2-n_{\sigma}^2x_u^1)
           \Big\} \\[0.5cm]
  &   & \negthickspace\displaystyle
        -\,\frac{1}{W^2}\,
           \sum_{\sigma=1}^{n-2}
           \Big\{L_{\sigma,11}^2+L_{\sigma,12}^2\Big\}\,(x_u^1x_v^2-x_u^2x_v^1) \\[0.6cm]
  &   & \negthickspace\displaystyle
        +\,\frac{1}{W}\,
           \sum_{\sigma=1}^{n-2}\sum_{\omega=1}^{n-2}
           \Big\{
             L_{\sigma,11}T_{\sigma,1}^\omega(n_\omega^1 x_v^2-n_\omega^2 x_v^1)
             -L_{\sigma,12}T_{\sigma,1}^\omega(n_\omega^1 x_u^2-n_\omega^2 x_u^1)
           \Big\} \\[0.5cm]
  &   & \negthickspace\displaystyle
        +\,\frac{1}{2W^2}\,
           \sum_{\sigma=1}^{n-2}
           \Big\{
             L_{\sigma,11}W_v(x_u^2n_\sigma^1-x_u^1n_\sigma^2)
             +L_{\sigma,11}W_u(x_v^2n_\sigma^1-x_v^1n_\sigma^2)\Big\} \\[0.5cm]
  &   & \negthickspace\displaystyle
        -\,\frac{1}{2W^2}\,
           \sum_{\sigma=1}^{n-2}
           \Big\{
             L_{\sigma,12}W_u(x_u^2n_\sigma^1-x_u^1n_\sigma^2)
             -L_{\sigma,12}W_v(x_v^2n_\sigma^1-x_v^1n_\sigma^2)
           \Big\} \\[0.5cm]
  &   & \negthickspace\displaystyle
        +\,\frac{1}{W}\,
           \sum_{\sigma=1}^{n-2}\sum_{\omega=1}^{n-2}
           \Big\{L_{\sigma,11}L_{\omega,12}-L_{\sigma,12}L_{\omega,11}\Big\}\,
           (n_\sigma^1n_\omega^2-n_\sigma^2n_\omega^1) \\[0.5cm]
  &   & \negthickspace\displaystyle
        +\,\frac{1}{W}\,
           \sum_{\sigma=1}^{n-2}
           \Big\{
             2H_{\sigma,u}W+2H_\sigma W_u
           \Big\}\,
           (n_\sigma^1x_v^2-n_\sigma^2x_v^1).
\end{array}
\end{equation}
An analog expression can be computed for $\chi_{vv}.$ Summation would show
\begin{equation}\label{3.8}
\begin{array}{lll}
  \triangle\chi\negthickspace
  & = & \negthickspace\displaystyle
        \frac{1}{W}\,
        \sum_{\sigma=1}^{n-2}
        \Big\{
          H_\sigma W_u(n_\sigma^1x_v^2-n_\sigma^2x_v^1)
          -H_\sigma W_v(n_\sigma^1x_u^2-n_\sigma^2x_u^1)
        \Big\} \\[0.6cm]
  &   & \negthickspace\displaystyle
        -\,\frac{1}{W^2}\,
           \sum_{\sigma=1}^{n-2}
           \Big\{
             L_{\sigma,11}^2+2L_{\sigma,12}^2+L_{\sigma,22}^2
           \Big\}\,(x_u^1x_v^2-x_u^2x_v^1) \\[0.6cm]
  &   & \negthickspace\displaystyle
        +\,\frac{1}{W}\,
           \sum_{\sigma=1}^{n-2}\sum_{\omega=1}^{n-2}
           \Big\{
             (L_{\sigma,11}-L_{\sigma,22})L_{\omega,12}
             -(L_{\omega,11}-L_{\omega,22})L_{\sigma,12}
           \Big\}\,(n_\sigma^1n_\omega^2-n_\sigma^2n_\omega^1) \\[0.6cm]
  &   & \negthickspace\displaystyle
        +\,\frac{1}{W}\,
           \sum_{\sigma=1}^{n-2}\sum_{\omega=1}^{n-2}
           \Big\{
             L_{\sigma,11}T_{\sigma,1}^\omega+L_{\sigma,12}T_{\sigma,2}^\omega
           \Big\}\,(n_\omega^1x_v^2-n_\omega^2x_v^1) \\[0.6cm]
  &   & \negthickspace\displaystyle
        +\,\frac{1}{W}\,
           \sum_{\sigma=1}^{n-2}\sum_{\omega=1}^{n-2}
           \Big\{
             -L_{\sigma,12}T_{\sigma,1}^\omega-L_{\sigma,22}T_{\sigma,2}^\omega
           \Big\}\,(n_\omega^1x_u^2-n_\omega^2x_u^1) \\[0.6cm]
  &   & \negthickspace\displaystyle
        +\,2\sum_{\sigma=1}^{n-2}
            \Big\{
              H_{\sigma,u}(n_\sigma^1x_v^2-n_\sigma^2x_v^1)
              -H_{\sigma,v}(n_\sigma^1x_u^2-n_\sigma^2x_u^1)
            \Big\}\,W \\[0.6cm]
  &   & \ldots\ldots\ldots\ldots\ldots\ldots\ldots\ldots\ldots\ldots\ldots\ldots\ldots
        \ldots\ldots\ldots\ldots\ldots\ldots\ldots\ldots\ldots\ldots
\end{array}
\end{equation}
  $$\begin{array}{l}
      \ldots\ldots\ldots\ldots\ldots\ldots\ldots\ldots\ldots\ldots\ldots
      \ldots\ldots\ldots\ldots \\[0.4cm]
      \displaystyle
      +\,\frac{1}{W}\,
         \sum_{\sigma=1}^{n-2}
         \Big\{
           (L_{\sigma,12,v}-L_{\sigma,22,u}
         \Big\}\,(n_\sigma^1x_v^2-n_\sigma^2x_v^1) \\[0.6cm]
      \displaystyle
      +\,\frac{1}{W}\,
         \sum_{\sigma=1}^{n-2}
         \Big\{
           L_{\sigma,11,v}-L_{\sigma,12,u}
         \Big\}\,(n_\sigma^1x_u^2-n_\sigma^2x_u^1).
    \end{array}\eqno(3.8)$$
In the last two rows we insert the equations of Codazzi/Mainardi (see \cite{Brauner_01}) in the form
\begin{equation}\label{3.9}
\begin{array}{lll}
  L_{\sigma,12,v}-L_{\sigma,22,u}\negthickspace
  & = & \negthickspace\displaystyle
        -\,H_\sigma W_u
        +\sum_{\omega=1}^{n-2}
         (L_{\omega,22}T_{\omega,1}^\sigma-L_{\omega,12}T_{\omega,2}^\sigma), \\[0.6cm]
  L_{\sigma,11,v}-L_{\sigma,12,u}\negthickspace
  & = & \negthickspace\displaystyle
        H_\sigma W_v
        +\sum_{\omega=1}^{n-2}
         (L_{\omega,12}T_{\omega,1}^\sigma-L_{\omega,11}T_{\omega,2}^\sigma),
\end{array}
\end{equation}
and we arrive at
\begin{equation}\label{3.10}
\begin{array}{lll}
  \triangle\chi\negthickspace
  & = & \negthickspace\displaystyle
        -\,\frac{1}{W^2}\,
           \sum_{\sigma=1}^{n-2}
           \Big\{
             L_{\sigma,11}^2+2L_{\sigma,12}^2+L_{\sigma,22}^2
           \Big\}\,
           (x_u^1x_v^2-x_u^2x_v^1) \\[0.6cm]
  &   & \negthickspace\displaystyle
        +\,\frac{1}{W}\,
           \sum_{\sigma=1}^{n-2}\sum_{\omega=1}^{n-2}
           \Big\{
             (L_{\sigma,11}-L_{\sigma,22})L_{\omega,12}
             -(L_{\omega,11}-L_{\omega,22})L_{\sigma,12}
           \Big\}\,
           (n_\sigma^1n_\omega^2-n_\sigma^2n_\omega^1) \\[0.6cm]
  &   & \negthickspace\displaystyle
        +\,2\sum_{\sigma=1}^{n-2}\sum_{\omega=1}^{n-2}
           \Big\{
             (n_\sigma^1x_v^2-n_\sigma^2x_v^1)T_{\sigma,1}^\omega
             -(n_\sigma^1x_u^2-n_\sigma^2x_u^1)T_{\sigma,2}^\omega
           \Big\}\,H_\sigma \\[0.6cm]
 &   & \negthickspace\displaystyle
        +\,2\sum_{\sigma=1}^{n-2}
            \Big\{
              H_{\sigma,u}(n_\sigma^1x_v^2-n_\sigma^2x_v^1)
              -H_{\sigma,v}(n_\sigma^1x_u^2-n_\sigma^2x_u^1)
            \Big\}\,W.
\end{array}
\end{equation}
The statement follows.
\end{itemize}
\end{proof}
\def\proofname{Proof}
\subsection{$\mu$-stability for mean curvature graphs}
Evaluating (\ref{3.2}) in general is rather difficult. Therefore, we confine to the following special situation:
\begin{theorem}
Given a conformally parametrized graph $X$ of prescribed mean curvature field $H=(H_1,\ldots,H_{n-2}),$ $H_\sigma=H_\sigma(X,Z),$ w.r.t. a torsion-free ONS ${\mathfrak N}.$ Then it holds
\begin{equation}\label{3.11}
\begin{array}{lll}
  \triangle\chi\negthickspace
  & = & \negthickspace\displaystyle
        -\,2(2H^2-K)W\chi
        +2\sum_{\sigma=1}^{n-2}
            (H_{\sigma,X}\cdot X_u^t+H_{\sigma,Z}\cdot N_{\sigma,u}^t)(n_\sigma^1x_v^2-n_\sigma^2x_v^1) \\[0.6cm]
  &   & \negthickspace\displaystyle
        -\,2\sum_{\sigma=1}^{n-2}
            (H_{\sigma,X}\cdot X_v^t+H_{\sigma,Z}\cdot N_{\sigma,v}^t)(n_\sigma^1x_u^2-n_\sigma^2x_u^1)
        \qquad\mbox{in}\ B.
\end{array}
\end{equation}
Assume, furthermore, that
\begin{equation}\label{3.12}
  0<\chi_{min}\le\chi\quad\mbox{in}\ B
\end{equation}
with a real constant $\chi_{min},$ and, furthermore, with real constants $h_{min},$ $h_1,$ and $h_2$ let
\begin{equation}\label{3.13}
  0<h_{min}\le H,\quad
  |H_{\sigma,X}|\le h_1\,,\quad
  |H_{\sigma,Z}|\le h_2
  \quad\mbox{for all}\ \sigma=1,\ldots,n-2.
\end{equation}
\goodbreak\noindent
Then, the graph is $\mu$-stable in the sense of
\begin{equation}\label{3.14}
  \int\hspace{-0.25cm}\int\limits_{\hspace{-0.3cm}B}|\nabla\varphi|\,dudv
  \ge\mu\int\hspace{-0.25cm}\int\limits_{\hspace{-0.3cm}B}(2H^2-K)W\varphi^2\,dudv
\end{equation}
for all $\varphi\in C_0^\infty(B,\mathbb R)$ and all
\begin{equation}\label{3.15}
  0<\mu\le 2\,-\frac{\sqrt{2}}{\chi_{min}}\,\left[\frac{(n-2)(h_1+h_2)}{h_{min}}+2h_2\right]
\end{equation}
with its mean curvature $H$ and its Gaussian curvature $K.$
\end{theorem}
\begin{remarks}\quad
\begin{itemize}
\item[1.]
If $X$ is a minimal surface, that is $h_{min},h_1,h_2=0,$ then condition (\ref{3.12}) is not needed. Note that because $\mu\le 2,$ all $\mu$-stabilities are included in the one stability $\mu=2.$
\item[2.]
For $n=3,$ condition (\ref{3.12}) means a bound on the gradient of the graph.
\end{itemize}
\end{remarks}
\def\proofname{Proof of the theorem}
\begin{proof}\quad
\begin{itemize}
\item[1.]
Because $X$ is a graph over the $[x,y]$-plane, the function $\chi$ is positive. We estimate as follows:
\begin{equation}\label{3.16}
\hspace*{-0.6cm}
\begin{array}{lll}
  \triangle\chi\negthickspace
  & \le & \negthickspace\displaystyle
          -\,2(2H^2-K)W\chi
          +\sqrt{32}\,\sum_{\sigma=1}^{n-2}h_1|X_u||X_v|
          +\sqrt{8}\,\sum_{\sigma=1}^{n-2}h_2(|N_{\sigma,u}||X_v|+|N_{\sigma,v}||X_u|) \\[0.6cm]
  & \le & \negthickspace\displaystyle
          -\,2(2H^2-K)W\chi
          +\sqrt{2}\,(n-2)(h_1+h_2)W
          +\sqrt{2}\,h_2\sum_{\sigma=1}^{n-2}|\nabla N_\sigma|^2 \\[0.6cm]
  & \le & \negthickspace\displaystyle
          -\,2(2H^2-K)W\chi
          +\sqrt{2}\,\left\{\frac{(n-2)(h_1+h_2)}{h_{min}}+2h_2\right\}(2H^2-K)W.
\end{array}
\end{equation}
\item[2.]
For any $\varphi\in C_0^\infty(B,\mathbb R)$ it holds $|\nabla\varphi|^2=\chi^2|\nabla(\varphi\chi^{-1})|^2+\frac{2}{\chi}\,\nabla\varphi\cdot\nabla\chi^t-\frac{\varphi^2}{\chi^2}\,|\nabla\chi|^2\,.$ Taking $\mbox{div}\left(\frac{\varphi^2\chi_u}{\chi}\,,\frac{\varphi^2\chi_v}{\chi}\right)=2\,\frac{\varphi}{\chi}\,\nabla\varphi\cdot\nabla\chi^t-\frac{\varphi^2}{\chi^2}\,|\nabla\chi|^2+\frac{\varphi^2}{\chi}\,\triangle\chi$ into account, we obtain
\begin{equation}\label{3.17}
\begin{array}{lll}
  |\nabla\varphi|^2\negthickspace
  &  =  & \negthickspace\displaystyle
          \chi^2|\nabla(\varphi\chi^{-1})|^2
          +\mbox{div}
           \left(
             \frac{\varphi^2\chi_u}{\chi}\,,\frac{\varphi^2\chi_v}{\chi}
           \right)
          -\frac{\varphi^2}{\chi}\,\triangle\chi \\[0.6cm]
  & \ge & \negthickspace\displaystyle
          \chi^2|\nabla(\varphi\chi^{-1})|^2
          +\mbox{div}
           \left(
             \frac{\varphi^2\chi_u}{\chi}\,,\frac{\varphi^2\chi_v}{\chi}
           \right)
          +2(2H^2-K)W\varphi^2 \\[0.6cm]
  &     & \negthickspace\displaystyle
          -\,\sqrt{2}\,
             \left\{
               \frac{(n-2)(h_1+h_2)}{h_{min}}+2h_2
             \right\}
             (2H^2-K)W\,\frac{\varphi^2}{\chi_{min}}\,.
\end{array}
\end{equation}
Integrating using conformal parameters shows
\begin{equation}\label{3.18}
\begin{array}{l}
  \displaystyle
  \int\hspace{-0.25cm}\int\limits_{\hspace{-0.3cm}B}
  \Big\{
    |\nabla\varphi|^2-\mu(2H^2-K)W\varphi^2
  \Big\}\,dudv \\[0.7cm]
  \hspace*{0.6cm}\displaystyle
  \ge\,\left\{
         2-\mu-\frac{\sqrt{2}}{\chi_{min}}\,\left[\frac{(n-2)(h_1+h_2)}{h_{min}}+2h_2\right]
       \right\}
       \int\hspace{-0.25cm}\int\limits_{\hspace{-0.3cm}B}(2H^2-K)W\,\varphi^2\,dudv,
\end{array}
\end{equation}
and the statement follows.
\end{itemize}
\end{proof}
\def\proofname{Proof}
\section{The Hopf field. Branch points}
\setcounter{equation}{0}
The Hopf-functions (cp. with \cite{HopfH_01}), which we now introduce, are essentially used to prove our $\mu$-stability result in the next section.
\begin{definition}
The Hopf function ${\mathcal H}_N$ of the surface along a unit normal vector $N$ is defined as
\begin{equation}\label{4.1}
  {\mathcal H}_N:=L_{N,11}-L_{N,22}-2iL_{N,12}\,.
\end{equation}
\end{definition}
\noindent
Recall Wirtinger's differential symbols $\Phi_w:=\frac{1}{2}\,(\Phi_u-i\Phi_v),$ $\Phi_{\overline w}:=\frac{1}{2}\,(\Phi_u+i\Phi_v)$ using complex variables $w=u+iv.$
\begin{proposition}
Given the section $(X,{\mathfrak N}),$ $X$ conformally parametrized, it holds
\begin{equation}\label{4.2}
  {\mathcal H}_{\sigma,\overline w}
  =2H_{\sigma,w}W
   +\sum_{\omega=1}^{n-2}
    \Big\{
      (L_{\omega,22}+iL_{\omega,12})T_{\omega,1}^\sigma
       -(L_{\omega,21}+iL_{\omega,11})T_{\omega,2}^\sigma
    \Big\}
\end{equation}
for all $\sigma=1,\ldots,n-2.$
\end{proposition}
\begin{remark}
If the ONS ${\mathfrak N}$ is free of torsion and $X$ is minimal, then all ${\mathcal H}_\sigma$ are holomorphic: $H_{\sigma,\overline w}\equiv 0.$
\end{remark}
\def\proofname{Proof of the Theorem}
\begin{proof}
For ${\mathcal H}_\sigma^*:=L_{\sigma,11}-iL_{\sigma,12}=\frac{1}{2}\,{\mathcal H}_\sigma+WH_\sigma$ it holds
\begin{equation}\label{4.3}
  {\mathcal H}_{\sigma,\overline w}^*
  =\frac{1}{2}\,(L_{\sigma,11,u}+L_{\sigma,12,v})
   +\frac{i}{2}\,(L_{\sigma,11,v}-L_{\sigma,12,u}).
\end{equation}
Inserting (\ref{3.9}) shows that
\begin{equation}\label{4.4}
\begin{array}{lll}
  {\mathcal H}_{\sigma,\overline w}^*\negthickspace
  & = & \negthickspace\displaystyle
        \frac{1}{2}\,{\mathcal H}_{\sigma,\overline w}
        +W_{\overline w}H_\sigma+WH_{\sigma,\overline w} \\[0.5cm]
  & = & \negthickspace\displaystyle
        \frac{1}{2}\,(W_u+iW_v)H_\sigma+WH_{\sigma,u} \\[0.5cm]
  &   & \negthickspace\displaystyle
        +\,\frac{1}{2}\,
           \sum_{\omega=1}^{n-2}
           \Big\{
             (L_{\omega,22}+iL_{\omega,12})T_{\omega,1}^\sigma
              -(L_{\omega,21}+iL_{\omega,11})T_{\omega,2}^\sigma
           \Big\},
\end{array}
\end{equation}
and rearranging the last identity for ${\mathcal H}_{\sigma,\overline w}$ proofs the statement.
\end{proof}
\def\proofname{Proof}
\begin{corollary}
Given the section $(X,{\mathfrak N})$ with a minimal surface $X$ together with a torsion-free ONS ${\mathfrak N}.$ Then, in every compact subset $B'\subset\subset B$ the Gaussian curvature of $X$ vanishes at most at isolated points.
\end{corollary}
\def\proofname{Proof of the Corollary}
\begin{proof}
There hold
\begin{equation}\label{4.5}
  |{\mathcal H}_\sigma|^2
  ={\mathcal H}_\sigma\cdot\overline{{\mathcal H}_\sigma}
  =4(-K_\sigma)W^2
  \quad\mbox{for all}\ \sigma=1,\ldots,n-2,
\end{equation}
and therefore
\begin{equation}\label{4.6}
  -K=\frac{1}{4W^2}\,\sum_{\sigma=1}^{n-2}|{\mathcal H}_\sigma|^2
\end{equation}
for the Gaussian curvature of the minimal immersion. But due to holomorphy, in every compact set $B'\subset\subset B$ all the $|{\mathcal H}_\sigma|,$ $\sigma=1,\ldots,n-2,$ vanish at most at isolated points.
\end{proof}
\def\proofname{Proof}
\section{$\mu$-stability due to Barbosa/do Carmo and Ruchert}
\setcounter{equation}{0}
\begin{theorem}
Given the section $(X,{\mathfrak N})$ with a conformally parametrized minimal immersion $X$ and with a torsion-free ONS ${\mathfrak N}.$ For any real $\kappa_0>0$ we assume that
\begin{equation}\label{5.1}
  Q:=\int\hspace{-0.25cm}\int\limits_{\hspace{-0.3cm}B}(\kappa_0-K)W\,dudv<\omega_0
\end{equation}
with a real constant $\omega_0\in(0,4\pi)$ and the Gaussian curvature $K$ of $X.$ Next, let $S_\omega^2\subset S^2$ denote a spherical cap such that
\begin{equation}\label{5.2}
  \mbox{\rm Area}\,(S_\omega^2)=\omega_0\,,
\end{equation}
and let $\mu>0$ be the smallest eigenvalue of the spherical Laplacian on $S_\omega^2.$ Then, $X$ is $\mu$-stable with this number $\mu,$ that is it holds
\begin{equation}\label{5.3}
  \int\hspace{-0.25cm}\int\limits_{\hspace{-0.3cm}B}|\nabla\varphi|^2\,dudv
  \ge\mu\int\hspace{-0.25cm}\int\limits_{\hspace{-0.3cm}B}(-K)W\varphi^2\,dudv
\end{equation}
for all $\varphi\in C_0^\infty(B,\mathbb R).$
\end{theorem}
\begin{proof}
With a real $\kappa_0>0$ we consider the auxiliary function $\gamma=\kappa_0-K>0$ (note $K\le 0$), and denote by $\widehat K$ the Gaussian curvature of the new regular and conformal metric with coefficients $\widehat g_{ij}:=\gamma g_{ij},$ $i,j=1,2.$ First, we will show $\widehat K\le 1$ using the methods of \cite{Ruchert_01}.
\begin{itemize}
\item[1.]
Due to the conformality it holds $K=-\frac{1}{W}\,\triangle\log\sqrt{W},$ such that
\begin{equation}\label{5.4}
  \gamma\widehat K
  =-\frac{1}{W}\,\triangle\log\sqrt{\gamma W}
  =K-\frac{1}{W}\,\triangle\log\sqrt{\gamma}
\end{equation}
and it follows
\begin{equation}\label{5.5}
  \gamma^3\widehat K=\gamma^2 K+\frac{2}{W}\,(\gamma_w\gamma_{\overline w}-\gamma\gamma_{w\overline w}).
\end{equation}
Analogously, we calculate
\begin{equation}\label{5.6}
  K=\frac{2}{W^2}\,\left\{\frac{W_wW_{\overline w}}{W}-W_{w\overline w}\right\}.
\end{equation}
\item[2.]
From (\ref{4.6}) we conclude
\begin{equation}\label{5.7}
  \frac{1}{4W^2}\,\sum_{\sigma=1}^{n-2}|{\mathcal H}_\sigma|^2
  =-K=\gamma-\kappa_0\,.
\end{equation}
\item[3.]
We introduce the following notation for an inner product: Let $v=(v_1,\ldots,v_{n-2})\in\mathbb C^{n-2}$ and $w=(w_1,\ldots,w_{n-2})\in\mathbb C^{n-2}$ be complex-valued vectors, $\overline v=(\overline v_1,\ldots,\overline v_{n-2})$ etc. for its complex conjugation, and
\begin{equation}\label{5.8}
  v\star w=w\star v=\sum_{\sigma=1}^{n-2}v_\sigma w_\sigma\in\mathbb C,\quad
  \|v\|^2=v\star\overline v
         =\sum_{\sigma=1}^{n-2}v_\sigma\overline v_\sigma=\sum_{\sigma=1}^{n-2}|v_\sigma|^2\in\mathbb R.
\end{equation}
Then we rewrite (\ref{5.7}) into the form
\begin{equation}\label{5.9}
  \gamma=\frac{1}{4W^2}\,{\mathcal H}\star\overline{\mathcal H}+\kappa_0
      =\frac{1}{4W^2}\,\|{\mathcal H}\|^2+\kappa_0\,.
\end{equation}
\item[4.]
We differentiate $\gamma$ (note that all ${\mathcal H}_{\sigma,\overline w}=0$ due to (\ref{4.2})) and obtain
\begin{equation}\label{5.10}
  \gamma_w
  =\frac{1}{4W^2}\,\overline{\mathcal H}\star{\mathcal H}_w
   -\frac{2W_w}{W}\,(\gamma-\kappa_0),\quad
  \gamma_{\overline w}
  =\frac{1}{4W^2}\,{\mathcal H}\star\overline{\mathcal H}_{\overline w}
   -\frac{2W_{\overline w}}{W}\,(\gamma-\kappa_0).
\end{equation}
Now the second derivatives:
\begin{equation}\label{5.11}
\begin{array}{lll}
  \gamma_{w\overline w}\negthickspace
  & = & \negthickspace\displaystyle
        \frac{\partial}{\partial\overline w}
        \left\{
          \frac{1}{4W^2}\,\overline{\mathcal H}\star{\mathcal H}_w
          -\frac{W_w}{2W^3}\,\|{\mathcal H}\|^2
        \right\} \\[0.6cm]
  & = & \negthickspace\displaystyle
        -\,\frac{1}{2W^3}\,
           \Big\{
             W_{\overline w}\,\overline{\mathcal H}\star{\mathcal H}_w
             +W_w\,{\mathcal H}\star\overline{\mathcal H}_{\overline w}
           \Big\}
        +\frac{1}{4W^2}\,\|{\mathcal H}_w\|^2 \\[0.6cm]
  &   & \negthickspace\displaystyle
        -\,\frac{1}{2W^3}
           \left\{
             W_{w\overline w}-\frac{3W_wW_{\overline w}}{W}
           \right\}\|{\mathcal H}\|^2\,.
\end{array}
\end{equation}
These identities imply
\begin{equation}\label{5.12}
\begin{array}{lll}
  \gamma_w\gamma_{\overline w}-\gamma\gamma_{w\overline w}\negthickspace
  & = & \negthickspace\displaystyle
        \frac{1}{16W^4}
        \left\{
         (\overline{\mathcal H}\star{\mathcal H}_w)
         ({\mathcal H}\star\overline{\mathcal H}_{\overline w})
         -\|{\mathcal H}\|^2\|{\mathcal H}_w\|^2
        \right\} \\[0.6cm]
  &   & \negthickspace\displaystyle
        -\,\frac{\|{\mathcal H}\|^2}{4W}\,K\gamma
        -\frac{\kappa_0}{W^2}\,\Big\|{\mathcal H}_w-\frac{2}{W}\,W_w{\mathcal H}\Big\|^2\,.
\end{array}
\end{equation}
\item[5.]
We estimate
\begin{equation}\label{5.13}
\begin{array}{lll}
  \displaystyle
  -\,\frac{\|{\mathcal H}\|^2}{4W}\,K\gamma\negthickspace
  &  =  & \negthickspace\displaystyle
          -\,(\gamma-\kappa_0)\gamma KW
          \,=\,-\,\frac{\gamma^2}{2}\,KW-\frac{\gamma^2}{2}\,KW+\kappa_0\gamma KW \\[0.5cm]
  &  =  & \negthickspace\displaystyle
          -\,\frac{\gamma^2}{2}\,KW
          -\frac{\gamma^2}{2}\,(\kappa_0-\gamma)W
          +\kappa_0\gamma(\kappa_0-\gamma)W \\[0.5cm]
  &  =  & \negthickspace\displaystyle
          \frac{\gamma^3}{2}\,W
          -\frac{\gamma^2}{2}\,KW
          +\gamma\Big(\kappa_0-\frac{3}{2}\,\gamma\Big)\kappa_0W
          \,\le\,\frac{\gamma^3}{2}\,W-\frac{\gamma^2}{2}\,KW,
\end{array}
\end{equation}
as well as
\begin{equation}\label{5.14}
  (\overline{\mathcal H}\star{\mathcal H}_w)({\mathcal H}\star\overline{\mathcal H}_{\overline w})
  -\|{\mathcal H}\|^2\|{\mathcal H}_w\|^2
  =(\overline{\mathcal H}\star{\mathcal H}_w)
   \overline{(\overline{\mathcal H}\star{\mathcal H}_w)}
   -\|{\mathcal H}\|^2\|{\mathcal H}_w\|^2
  \le 0
\end{equation}
using the Cauchy-Schwarz inequality.
\item[6.]
Finally, we arrive at $\gamma_w\gamma_{\overline w}-\gamma\gamma_{w\overline w}\le\frac{\gamma^3}{2}\,W-\frac{\gamma^2}{2}\,KW,$ therefore
\begin{equation}\label{5.15}
  \gamma^3\widehat K\le\gamma^2K+\gamma^3-\gamma^2K=\gamma^3\,,
\end{equation}
which means $\widehat K\le 1.$
\item[7.]
Now, we prove the $\mu$-stability inequality. Let $\widehat\triangle$ denote the Laplacian w.r.t. the regular metric $\widehat g_{ij},$ and let $\widehat\lambda>0$ be the first eigenvalue of the homogeneous Dirichlet problem
\begin{equation}\label{5.16}
  \widehat\triangle\varphi+\lambda\varphi=0\quad\mbox{in}\ B,\quad
  \varphi=0\quad\mbox{on}\ \partial B.
\end{equation}
Next, with $S_{\widetilde\omega}^2\subset S^2,$ $\widetilde\omega\in(0,4\pi),$ we mean a spherical cap with the property
\begin{equation}\label{5.17}
  \mbox{Area}\,(S_{\widetilde\omega}^2)=Q
\end{equation}
with $Q<\omega_0$ given in the Theorem, and $\lambda_1^*>0$ is the first eigenvalue of the problem
\begin{equation}\label{5.18}
  \triangle^*\varphi^*+\lambda^*\varphi^*=0\quad\mbox{in}\ S_{\widetilde\omega}^2\,,\quad
  \varphi^*=0\quad\mbox{on}\ \partial S_{\widetilde\omega}^2
\end{equation}
with the spherical Laplacian $\triangle^*.$
\item[8.]
Note that $\widehat K\le 1=K(S^2)$ for the Gaussian curvature of the sphere $S^2,$ and $Q$ is also the area w.r.t. the metric $\widehat g_{ij}.$ Because $S_{\widetilde\omega}^2$ is a geodesic disc in $S^2,$ it minimizes the first eigenvalue (see \cite{Barbosa_doCarmo_02}) such that
\begin{equation}\label{5.19}
  \lambda_1^*\le\widehat\lambda_1\,.
\end{equation}
On the other hand, by assumption, $S_\omega^2\subset S^2$ is a spherical cap with first eigenvalue $\mu$ and $\mbox{Area}\,(S_\omega^2)=\omega_0,$ and by $Q<\omega_0$ and the monotonicity of the first eigenvalue we have $\mu<\lambda_1^*$ (see \cite{Barbosa_doCarmo_02}). Using conformal parameters, this gives
\begin{equation}\label{5.20}
  \mu
  <\lambda_1^*
  \le\widehat\lambda_1
  \le\frac{\displaystyle\int\hspace{-0.25cm}\int\limits_{\hspace{-0.3cm}B}|\nabla\varphi|^2\,dudv}
          {\displaystyle\int\hspace{-0.25cm}\int\limits_{\hspace{-0.3cm}B}(\kappa_0-K)W\varphi^2\,dudv}
  \quad\mbox{for all}\ \varphi\in C_0^\infty(B,\mathbb R),
\end{equation}
while rearranging shows
\begin{equation}\label{5.21}
  \int\hspace{-0.25cm}\int\limits_{\hspace{-0.3cm}B}|\nabla\varphi|^2\,dudv
  \ge\mu\int\hspace{-0.25cm}\int\limits_{\hspace{-0.3cm}B}(\kappa_0-K)W\varphi^2\,dudv
  \ge\mu\int\hspace{-0.25cm}\int\limits_{\hspace{-0.3cm}B}(-K)W\varphi^2\,dudv
\end{equation}
for all $\varphi\in C_0^\infty(B,\mathbb R)$ due to $\kappa_0>0$ and $K\le 0.$ The statement is proved.
\end{itemize}
\end{proof}
\begin{remark}
With a real $a\in(0,2],$ let $\widehat K\le a$ for the Gaussian curvature $\widehat K$ in our proof. Due to \cite{Barbosa_doCarmo_02}, a minimal surface is stable ($\mu=2$) if
\begin{equation}\label{5.22}
  \int\hspace{-0.25cm}\int\limits_{\hspace{-0.3cm}B}(-K)W\,dudv\le\frac{4\pi}{1+a}\,.
\end{equation}
In fact, $\widehat K\le 2$ is established in \cite{Barbosa_doCarmo_02} with different methods than we used above. We could improve this curvature estimate ($\widehat K\le 1$) under the assumption of a flat normal bundle.
\end{remark}

\vspace*{0.8cm}
\noindent
Steffen Fr\"ohlich\\
Freie Universit\"at Berlin\\
Fachbereich Mathematik und Informatik\\
Institut f\"ur Mathematik I\\
Arnimalle 2-6\\
D-14195 Berlin\\
Germany\\[0.2cm]
e-mail: sfroehli@mi.fu-berlin.de

\end{document}